\documentclass{elsart3-1}

 \usepackage{graphicx}

\usepackage{amssymb}

\usepackage[english,francais]{babel}
\usepackage[T1]{fontenc}
\usepackage[latin1]{inputenc}

\newtheorem{e-proposition}[theorem]{Proposition}

\newtheorem{e-definition}[theorem]{Definition\rm}

\newtheorem{theoreme}{Th\'eor\`eme}[section]
\newtheorem{lemme}[theoreme]{Lemme}

\newtheorem{definition}[theoreme]{D\'efinition\rm}

\renewcommand{\theequation}{\arabic{equation}}
\def\og{\leavevmode\raise.3ex\hbox{$\scriptscriptstyle\langle\!\langle$~}}
\def\fg{\leavevmode\raise.3ex\hbox{~$\!\scriptscriptstyle\,\rangle\!\rangle$}}

\begin{document}
\centerline{Statistique}
\begin{frontmatter}




%
\selectlanguage{francais}
\title{Consistance d'un estimateur de minimum de variance \'etendue}

\vspace{-2.6cm}
\selectlanguage{english}
\title{Consistency of a least extended variance estimator}



\author[authorlabel1]{Joseph Rynkiewicz}
\ead{joseph.rynkiewicz@univ-paris1.fr}
\address[authorlabel1]{SAMOS/MATISSE, Universit\'e de ParisI, 72 rue Regnault 75013 Paris, France, Tél. et Fax : 01-44-07-87-05}

\begin{abstract}
We consider a generalization of the criterion minimized by the K-means algorithm, where a neighborhood structure is used in the calculus of the variance. 
Such tool is used, for example with Kohonen maps, to measure the quality of the quantification preserving the neighborhood relationships.
If we assume that the parameter vector is in a compact Euclidean space and all it components are separated by a minimal distance,  we show the strong consistency of the set of parameters almost realizing  the minimum of the empirical extended variance. {\it To cite this article:}

\vskip 0.5\baselineskip

\selectlanguage{francais}
\noindent{\bf R\'esum\'e}
\vskip 0.5\baselineskip
\noindent
On consid\`ere une g\'en\'eralisation du critère minimisé par l'algorithme des K-moyennes [K-means], où une structure de voisinage est introduit dans le calcul de la variance.  Un tel outil est utilis\'e, par exemple avec des cartes de Kohonen,  pour mesurer la qualité de la quantification respectant les structures de voisinage. Si on suppose que le vecteur param\`etre est dans un compact d'un espace euclidien et que toutes ses composantes sont s\'epar\'ees par une distance minimale, on montre la consistance forte de l'ensemble des param\'etres assez proches du minimum de variance \'etendue. {\it Pour citer cet article~:}

\end{abstract}
\end{frontmatter}

\selectlanguage{francais}
\section{Introduction}
Nous consid\'erons une g\'en\'eralisation de la variance intra-classe qui est considérée comme le principal critère de mesure de qualité des cartes de Kohonen (cf Kohonen \cite{Kohonen}), bien que l'algorithme de Kohonen ne minimise pas exactement ce critère (cf Cottrell et al. \cite{Cottrell}). La variance étendue est la somme de la variance intra-classe et d'un terme qui d\'epend des classes voisines. Sa minimisation permet notamment d'obtenir une classification qui respecte les relations de voisinage et qui donne lieu \`a des interpr\'etations ais\'ees, puisque la proximit\'e des classes correspond \`a la proximit\'e des données r\'eelles représentées par ces classes.

Nous considérons dans toute la suite que les observations i.i.d. sont dans le compact  $[0,1]^d$, qu'elles ont pour mesure de probabilité $P$ qui admet une densité par rapport à la mesure de Lebesque, bornée par une constante $B$. 
Dans la suite, on appelera ``centro\"ide'' un vecteur de $\left[0,1\right]^d$ qui représente une classe d'observations $\omega$.
\begin{definition}  
Pour $e\in\mathbb N^*, e\leq d$, soit un ensemble fini $I\subset \mathbb Z^e$ et $\Lambda$ la fonction de voisinage d\'efinie de $I-I := \{i-j, i, j \in I\}$ dans $\left[0, 1\right]$ telle que $\Lambda(k)=\Lambda(-k)$ et $\Lambda(0)=1$.
\end{definition} 
\begin{definition} 

Notons $\Vert .\Vert $ la norme euclidienne, soit 
\[
D_{I}^{\delta}:=\left\{x:=\left(x_i\right)_{i \in I} \in \left(\left[0,1\right]^d\right)^I \mbox{ tels que } \left \Vert x_i-x_j\right \Vert \geq \delta, \mbox{ si } i\neq j\right\}
\]
l'ensemble des centroïdes $x_i$ séparés par une distance d'au moins $\delta$.
\end{definition}
\begin{definition}
La tessellation de Voronoï $\left(C_i(x)\right)_{i\in I}$ est définie par 
\[
C_i(x):=\left\{\omega \in \left[0,1\right]^d  \mbox{ tels que }\Vert x_i-\omega \Vert < \Vert x_j - \omega \Vert \mbox{ si } j \neq i \right\}
\]
En cas d'ex-aequo, on assigne $\omega \in C_i(x)$ grâce à l'ordre lexicographique sur $I$.
On remarquera alors que $\left(C_i(x)\right)_{i\in I}$ est une partition borélienne  dont l'intérieur coincide avec la mosa\"ique de Vorono\"i ouverte. Réciproquement, l'indice de la tesselation de Voronoï pour une observation $\omega$ est définie par
\[
C_{x}^{-1}(\omega):=i \in I \mbox{ tel que }\omega \in C_i(x) 
\]
\end{definition}

\begin{definition}
La variance étendue est : 
\(
V(x):=\frac{1}{2}\sum_{i,j\in I}\Lambda(i-j)\int_{C_i(x)}\Vert x_j-\omega \Vert^2 dP(\omega)
\)\\  
De même, lorsqu'il y a un nombre fini $n$  d'observations, on définit la variance étendue empirique :\\
\(
V_n(x):=\frac{1}{2n}\sum_{i\in I}\sum_{\omega \in C_i(x)}\left(\sum_{j \in I}\Lambda(i-j)\Vert x_j - \omega \Vert^2\right)
\)
\end{definition}
Si une observation se trouve sur un hyperplan médiateur entre deux centro\"ides, tout déplacement d'un de ces deux centroïdes entrainera un saut de la variance étendue à moins que la mesure de cette observation ne soit nulle (par exemple, si la mesure ne charge pas les traces d'hyperplan). La fonction de variance étendue empirique $V_n(x)$ n'est donc pas continue et il n'existe pas, en général, d'ensemble de centro\"ides réalisant son minimum. Cependant, si on considère les suites $x^n$ telles que $V_n(x^n)$ soit suffisamment proche de son minimum, on peut se demander si ces suites convergent vers l'ensemble des centroïdes minimisant la variance théorique $V(x)$. Pour cela, nous procédons selon le même schéma de démonstration que Pollard \cite{Pollard} et nous commençons par montrer que les fonctions de variance étendue  vérifient une loi uniforme des grands nombres.
\section{Loi uniforme des grands nombres}  

Soit la famille de fonctions 
\[
{\mathcal G}:=\left\{ g_x(\omega):=\sum_{j\in I} \Lambda\left(C_x^{-1}(\omega)-j\right)\Vert x_j-\omega\Vert^2\mbox{ pour } x\in D_I^\delta\right\}
\]

Pour montrer la loi uniforme des grands nombres, il suffit de montrer que 
\begin{equation}
\sup_{x\in D^\delta_I}\left|\int g_x(\omega)dP_n(\omega)-\int g_x(\omega)dP(\omega)\right|\stackrel{p.s.}{\longrightarrow}0
\label{lgnu}\end{equation}
puisque, pour toute mesure de probabilité $Q$ sur $[0,1]^d$:
\[
\int g_x(\omega)dQ(\omega)=\int \sum_{j\in I}\Lambda\left(C_x^{-1}(\omega)-j\right)\Vert x_j-\omega \Vert^2dQ(\omega)=
\frac{1}{2}\sum_{i,j\in I}\Lambda(i-j)\int_{C_i(x)}\Vert x_j-\omega \Vert^2 dQ(\omega)
\]
d'après Gaenssler et Stute \cite{Gaenssler}, une condition suffisante pour que l'équation (\ref{lgnu}) soit vérifiée est que :
$\forall \varepsilon > 0, \forall x_0 \in D_I^\delta$ il existe un voisinage $S(x_0)$ de $x_0$ tel que 
\[
\int g_{x_0}(\omega)dP(\omega)-\varepsilon<\int\left(\inf_{x\in S\left(x_0\right)}g_x(\omega)\right)dP(\omega)
\leq \int \left( \sup_{x\in S(x_0)} g_x(\omega)\right)dP(\omega)<\int g_{x_0}(\omega)dP(\omega)+\varepsilon
\]
On peut d'abord prouver le résultat suivant, en utilisant une technique similaire à la preuve du lemme 11 de Fort et Pagès \cite{Fort}. 

\begin{lemme}\label{lemme1}
Soit $x\in D^\delta_I$ et $\lambda$ la mesure de Lebesgue sur $[0,1]^d$. Notons $E^c$ le complémentaire de l'ensemble $E$ dans $[0,1]^d$ et $|I|$ le cardinal de l'ensemble I. Pour $0<\alpha<\frac{\delta}{2}$, soit\\
\(
U^\alpha_i(x)=\left\{\omega \in [0,1]^d/\exists y\in D^\delta_I, x_j=y_j\mbox{ si } j\neq i\mbox{ et }\Vert x_i-y_i\Vert<\alpha \mbox{ et } \omega \in C_i^c(y)\cap C_i(x)\right\}
\)\\
l'ensemble des $\omega$ changeant de cellule de Voronoï lorsque le centroïde $x_i$ se déplace d'une distance d'au plus $\alpha$. Alors
\[
sup_{x\in D_I^\delta}\lambda\left(U^\alpha_i(x)\right)< \left(|I|-1\right)\left(\frac{2\alpha}{\delta}+\alpha\right)\left(\sqrt{2}\right)^{d-1}
\]

\end{lemme}

Considérons maintenant $x^0\in D_I^\delta$ et $S(x^0)$ un voisinage de $x^0$ inclus dans une boule de rayon $\alpha$, pour la distance euclidienne sur $D^\delta_I$. Soit $W(x^0)$ l'ensemble des $\omega$ restant dans leur cellule de Voronoï lorsque on déplace  $x^0$ vers n'importe quel $x\in S(x_0)$. Pour tout $\omega \in W(x^0)$ on a
\[
\begin{array}{l}
\inf_{x\in S(x^0)} g_x(\omega)\geq g_{x^0}(\omega)-\sum_{j\in I}\Lambda \left(C_{x^0}^{-1}(\omega)-j\right)\left(\Vert x^0_j-\omega \Vert^2 - \inf_{x\in S(x^0)} \Vert x^0_j-\omega \Vert^2 \right)\\
\geq g_{x_j^0}(\omega)-\sum_{j\in I}\left(\Vert x_j^0-\omega \Vert^2 - \inf_{x\in S(x^0)} \Vert x_j^0-\omega \Vert^2 \right)
\end{array}
\]  
Pour tout $\omega\in[0,1]^d$, on a, pour $\alpha$ suffisamment petit, 
\(
\left( \Vert x_j^0-\omega \Vert^2-\inf_{x\in S\left(x^0\right)}\Vert x_j-\omega\Vert^2\right)<\frac{\varepsilon}{2B|I|}
\) 
ainsi
\[
\int_{W(x^0)}\sum_{j\in I}\left(\Vert x_j^0-\omega\Vert^2-\inf_{x\in S\left(x^0\right)}\Vert x_j-\omega\Vert^2\right)dP(\omega)<\frac{\varepsilon}{2}
\mbox{ et }
\int_{W(x^0)}\left( g_{x^0}(\omega)-\inf_{x\in S(x^0)}g_x(\omega) \right)<\frac{\varepsilon}{2}
\]

Soit, maintenant $W(x^0)^c$, l'ensemble des $\omega$ changeant de cellule de Voronoï quand les centroïdes vont de $x^0$ vers un $x\in S_{x^0}$.  Si $\alpha<\frac{\delta}{2|I|}$, alors en déplaçant séquentiellement les composantes $x^0_i$ de $x^0$ vers $x_i$ de $x$, chaque configuration intermédaire reste dans $D^{\frac{\delta}{2}}_I$. Comme, pour tout $i\in I$, $\Vert x_i-\omega\Vert^2$ est borné par $1$ sur $[0,1]^d$, le lemme \ref{lemme1}, assure alors que

\[
\int_{W(x^0)^c}g_x(\omega)dP(\omega)<B|I|(|I|-1))\left(\frac{4\alpha}{\delta}+\alpha\right)\left(\sqrt{2}\right)^{d-1}
\]
Finalement, si on choisit $\alpha$ suffisamment petit pour que $B|I|(|I|-1))\left(\frac{4\alpha}{\delta}+\alpha\right)\left(\sqrt{2}\right)^{d-1}<\frac{\varepsilon}{2}$, on obtient
\[
\int_{D_I^\delta}g_{x^0}(\omega)dP(\omega)-\varepsilon<\int_{D_I^\delta}\left(\inf_{x\in S(x^0)}g_x(\omega)\right)dP(\omega)
\]
Exactement de la même façon, pour $\alpha$ suffisamment petit, on obtient :
\[
\int_{D_I^\delta}\left(\sup_{x\in S(x^0)}g_x(\omega)\right)dP(\omega)<\int_{D_I^\delta}g_{x^0}(\omega)dP(\omega)+\varepsilon
\]
Ainsi, la condition suffisante pour la loi uniforme des grands nombres est vraie pour la variance étendue.
\section{Consistance}
On veut montrer la consistance des centro\"ides qui minimisent ``presque'' la variance étendue dans $D_I^\delta$. 
Soit  l'ensemble des ``quasi-estimateurs'' de minimum de variance étendue : 
\[
\bar{\chi}_n^\beta:=\left\{x\in D_I^\delta\mbox{ tels que } V_n(x)<\inf_{x\in D_I^\delta}V_n(x)+\frac{1}{\beta(n)}\right\}
\] 
avec $\beta(n)$ une fonction strictement positive tel que $\lim_{n\rightarrow+\infty}\beta(n)=\infty$.
Soit
\(
\bar{\chi}=\arg\min_{x\in D_I^\delta}V(x)
\) l'ensemble qui minimise la variance étendue théorique, comme la fonction 
\(
 x\longmapsto V\left(x\right)
\)
est continue et non constante sur $D_{I}^{\delta }$, pour tout voisinage
$\mathcal{N}$ de $\bar{\chi }$, il existe $\eta \left(\mathcal{N}\right)>0$
tel que 
\[
\forall x\in D_{I}^{\delta }\backslash \mathcal{N},\, V\left(x\right)>\min _{x\in D_{I}^{\delta }}V\left(x\right)+\eta \left(\mathcal{N}\right)\]
Pour montrer la consistance forte, il suffit de montrer que pour tout voisinage
$\mathcal{N}$ de $\bar{\chi }$ on a \[
\lim _{n\rightarrow \infty }\bar{\chi }_{n}^{\beta }\stackrel{p.s}{\subset }\mathcal{N}\Longleftrightarrow \lim _{n\rightarrow \infty }V\left(\bar{\chi }_{n}^{\beta }\right)-V\left(\bar{\chi }\right)\stackrel{p.s.}{\leq }\eta \left(\mathcal{N}\right)\]
avec  
\(
V\left(E\right)-V\left(F\right):=\sup \left\{ V\left(x\right)-V\left(y\right)\mbox {\, pour\, }x\in E\mbox {\, et\, }y\in F\right\} 
\).

Par définition $V_n\left(\bar{\chi }_{n}^{\beta }\right)\stackrel{p.s.}{\leq} V_{n}\left(\bar{\chi }\right)+\frac{1}{\beta \left(n\right)}$, 
de plus la loi uniforme des grands nombres assure que\\
\(
\lim _{n\rightarrow \infty }V_{n}\left(\bar{\chi }\right)-V\left(\bar{\chi }\right)\stackrel{p.s}{=}0
\), on obtient ainsi \(
\lim _{n\rightarrow \infty }V_{n}\left(\bar{\chi }_{n}^{\beta }\right)\stackrel{p.s.}{\leq }V\left(\bar{\chi }\right)+\frac{\eta \left(\mathcal{N}\right)}{2}\), 
de même on aura  $\lim _{n\rightarrow \infty }V\left(\bar{\chi }_{n}^{\beta }\right)-V_{n}\left(\bar{\chi }_{n}^{\beta }\right)\stackrel{p.s.}{=}0$
 et
\[
\lim _{n\rightarrow \infty }V\left(\bar{\chi }_{n}^{\beta }\right)-\frac{\eta \left(\mathcal{N}\right)}{2}\stackrel{p.s.}{<}\lim _{n\rightarrow \infty }V_{n}\left(\bar{\chi }_{n}^{\beta }\right)\stackrel{p.s.}{\leq }V\left(\bar{\chi }\right)+\frac{\eta \left(\mathcal{N}\right)}{2}\]
finalement 
\(
\lim _{n\rightarrow \infty }V\left(\bar{\chi }_{n}^{\beta }\right)-V\left(\bar{\chi }\right)\stackrel{p.s.}{\leq }\eta \left(\mathcal{N}\right)
\)
ce qui prouve la consistance forte du quasi-estimateur de minimum de variance étendue. 

\end{document}